\documentclass[12pt]{article}
\usepackage{amsmath,amssymb,amsthm}
\usepackage{latexsym}

\newcommand{\bysame}{%
\leavevmode\hbox to 3em{\hrulefill}\,}
\newcommand{\sect}[1]{\vspace{1cm}\noindent {\bf\large {#1}} \vspace{.5cm}}
\newcommand{\subsect}[1]{\noindent {\bf {#1}} \,\,}

\newcommand{\Z}{\mathbb{Z}}

\newcommand{\C}{\mathbb{C}}
\newcommand{\F}{\mathbb{F}}

\begin{document}

\begin{center}
{\bf\large Gauss Sums
 on $GL_2(\Z/p^l\Z)$}\\
\vspace{0.5cm}
{\large Taiki Maeda}\footnote{\it E-mail address: b0775003@sophia.jp}\\
{\it{\footnotesize Department of Mathematics, Sophia University, Kioicho 7-1, Chiyodaku, Tokyo 102-8554, Japan}}
\end{center}

\begin{abstract}
We determine explicitly the Gauss sum $\tau_l(\chi,e)=\sum_{X \in GL_2(\Z/p^l\Z)} \chi(X) e(\mathrm{Tr}X)$ 
on the general linear group $GL_2(\Z/p^l\Z)$ for every irreducible character $\chi$ of $GL_2(\Z/p^l\Z)$ and a nontrivial additive character $e$ of $\Z/p^l\Z$, 
where $p$ is an odd prime and $l$ is an integer $\geqslant 2$. While there are several studies of the Gauss sums on finite algebraic groups defined over a finite field, 
this paper seems to be the first one which determines the Gauss sums on a matrix group over a finite ring.
\end{abstract}

\sect{1. Introduction}

Generalizations of Gauss sums have been considered in many ways. 

As for the Gauss sums on finite algebraic groups, 
Kondo \cite{k} firstly determined the value of Gauss sum on the finite general linear group $GL_n(\F_q)$ for every irreducible character, 
where $\F_q$ denotes a finite field with $q$ elements. Also in a series of papers starting with \cite{kl}, Kim-Lee, D. S. Kim and Kim-Park 
described the values of Gauss sums on finite classical groups for linear characters. Saito-Shinoda \cite{ss1, ss2} considered the Gauss sums on
finite reductive groups for the Deligne-Lusztig generalized characters and applied this result to determine the value of Gauss sum on the finite symplectic group 
$Sp(4,q)$ and considered Gauss sums on the Chevalley group of type $G_2$ corresponding to every unipotent character. 

For an arbitrary finite group $G$, Gomi-Maeda-Shinoda \cite{gms} defined the Gauss sums on $G$, and explicitly determined the values of Gauss sums 
on the complex reflection group $G(m,r,n)$ and also on the Weyl group of any type for all irreducible characters. 
 
On the other hand, Gauss sums over $\Z/p^l\Z$ ($p$:odd prime, $l \geqslant 2$) were determined explicitly by Odoni \cite{od} and also by Funakura \cite{f} for $p=2$. 
To be precise, the Gauss sum over $\Z/p^l\Z$ is given as follows:
\[
g_l(\mu, e)=\sum_{x\in(\Z/p^l\Z)^{\times}}\mu (x)e(x) ,
\]
where $\mu$ is a multiplicative character of $(\Z/p^l\Z)^{\times}$ and $e$ is a nontrivial additive character of $\Z/p^l\Z$. 

Thus, in these vein, it would be quite natural to consider the Gauss sum on $G(\Z/p^l\Z)$, where $G$ is an algebraic group. 
To start to consider this problem, we deal with the case $G=GL_2$ in this paper. 
Extending the definition of the sum $g_l(\mu, e)$, we define the Gauss sum on the general linear group $GL_2(\Z/p^l\Z)$ as follows:
\[
\tau_l(\chi,e)= \sum_{X \in GL_2(\Z/p^l\Z)}\chi(X) e(\mathrm{Tr}X),  
\]
where $\chi$ is a character of $GL_2(\Z/p^l\Z)$. 
The definition of this Gauss sum differs from that of Gauss sum in \cite{gms} in the sense that $e$ is an additive character of $\Z/p^l\Z$. 
We will determine explicitly $\tau_l(\chi, e)$ for every irreducible character $\chi$ in the case that $p$ is an odd prime.  

In order to evaluate $\tau_l(\chi, e)$, we need the character theory of $GL_2(\Z/p^l\Z)$. 
The irreducible characters of $GL_2(\Z/p^l\Z)$ were completely constructed by Nobs \cite{n} and Leigh-Cliff-Wen \cite{lcw}. 
Nobs constructed them by using the Weil representation. However he did not consider the problem of finding the character values. 
On the other hand, by making full use of Clifford theory, Leigh-Cliff-Wen determined the values of the irreducible characters of $GL_2(\Z/p^l\Z)$. 
They also showed that the degrees of irreducible characters of $GL_2(\Z/p^l\Z)$ which do not come from $GL_2(\Z/p^{l-1}\Z)$ 
have precisely three possibilities, namely $p^{l-1}(p+1)$, $p^{l-1}(p-1)$, and $p^{l-2}(p^2-1)$, and such irreducible characters are induced from 
a character $\psi$ of the stabilizer of a linear character defined on a congruence subgroup $K_n$. 
However their construction of $\psi$ in some cases are not sufficient enough for us to compute the Gauss sum. 
Hence, for those cases, we give a more explicit construction and in particular, in Lemma 3.4, we construct $\psi$ in a different way. 

This paper is organized as follows:
in $\S2$, we recall the construction of the characters of $GL_2(\Z/p^l\Z)$ after \cite{lcw}. 
In $\S3$, we explicitly determine the Gauss sums on $GL_2(\Z/p^l\Z)$ for all irreducible characters after constructing the $\psi$'s. 
In $\S4$, we clarify some vague values of 
$P=\sum_{0\leqslant c,d <p^{i}, p\nmid d}\lambda (p^j\beta d+d^{-1}(p^kb-c^2))$ in (\cite{lcw}, 6.1), since we heavily use the results in \cite{lcw}. 

\sect{2. Preliminaries}

\subsect{2.1.} We shall use the following usual notation in the character theory of finite groups. 
Let $G$ be a finite group. The set of complex irreducible characters of $G$ is denoted by $\widehat{G}$. 
For $\C$-valued functions $f$ and $g$ on $G$, let
\[
\langle f, g \rangle_{G} = \frac{1}{|G|}\sum_{x\in G} f(x)\overline{g(x)}
\]
be the usual hermitian inner product on the vector space of $\C$-valued functions on $G$, where $\overline{g(x)}$ is the complex conjugate of $g(x)$.

Now we recall Clifford theory (cf. \cite{isa}). 
Let $N$ be a normal subgroup of $G$. For a character $\phi$ of $N$ and $g\in G$, we define $\phi^g$, $\phi^g:N \to \C$, by
$\phi^g(n)=\phi(gng^{-1})$. Then $\phi^{g}$ is also a character of $N$. 
Let $T=\mathrm{Stab}_{G}\phi=\{g \in G\ | \ \phi^g= \phi \}$ be the stabilizer of $\phi$ in $G$.
If $\phi\in\widehat{N}$, then Clifford theory implies that the map $\psi\mapsto \mathrm{ind}_{T}^{G}\psi$ is a bijection from
$\{\psi \in \widehat{T}\ | \ \langle \psi|_{N}, \phi \rangle_{N}>0\}$ onto $\{\chi \in \widehat{G}\ | \ \langle \chi|_{N}, \phi \rangle_{N}>0\}$. \\
 
\subsect{2.2.} We recall the construction of the characters of $GL_2(\Z/p^l\Z)$, following \cite{lcw}.
Let $p$ be an odd prime, $l\geqslant 2$ a positive integer, $m=\lfloor l/2\rfloor$, $n=\lceil l/2 \rceil$, 
$\Z_{p^l}=\Z/p^l\Z$ with its multiplicative group $\Z_{p^l}^{\times}$, and $G_l=GL_2(\Z_{p^l})$. We denote by $I$ the identity matrix of $G_l$ 
and let $K_{i}=\{I+p^iB \ | \ B \in M_2(\Z_{p^l}) \}$ for $1 \leqslant i< l$. 
Note that $K_i$ is a normal subgroup of $G_l$ for all $i$, and that $K_i$ is abelian if $i \geqslant n$. 
For a positive integer $k$, we will write $\zeta_k=\exp(2{\pi}\sqrt{-1}/k)$.
Let $\lambda$ be an injective additive character of $\Z_{p^l}$ defined by $\lambda(1)=\zeta_{p^l}$.
For $A\in M_2(\Z_{p^l})$, we define a character $\phi_A$ of $K_{n}$, $\phi_A:K_n \to \C^{\times}$, by 
\[
\phi_A(X)=\lambda(\mathrm{Tr}(A(X-I))).  
\]
For $\alpha \in \Z_{p^l}$, let $\mu_\alpha$ be a multiplicative character of $\Z_{p^l}^{\times}$, $\mu_{\alpha}:\Z_{p^l}^{\times}\to \C^{\times}$, such that
\[
\mu_\alpha (1+p^n)=\lambda (p^n\alpha ).
\]

Following \cite{lcw}, let  
\[
\begin{split}
\mathfrak{X}&=\{\chi \in\widehat{G_l}\ |\ K_{l-1}\nsubseteq \ker\nu\chi\ \textrm{for all}\ \nu \in \mathrm{Hom}(G_l, \C^{\times})  \}, \\
\mathfrak{X}_1&=\{ \chi \in \widehat{G_l}\ |\ \chi(I)=p^{l-1}(p+1)\}, \\
\mathfrak{X}_2&=\{ \chi \in \widehat{G_l}\ |\ \chi(I)=p^{l-1}(p-1)\}, \\
\mathfrak{X}_3&=\{ \chi \in \widehat{G_l}\ |\ \chi(I)=p^{l-2}(p^2-1)\}, \\
\mathfrak{X}_4&=\widehat{G_l}-\mathfrak{X}. 
\end{split}
\]
Also for $A\in M_2(\Z_{p^l})$, define
\[
\widehat{G_l}(A)=\{ \chi \in \widehat{G_l}\ |\ \langle \chi |_{K_{n}}, \phi_A\rangle_{K_{n}} >0 \}. 
\]

The following properties were proven in \cite{lcw} (Theorem 3.1, its proof, and 3.5).\\

\subsect{Theorem 2.1.}
\textit{With the notation as above, we have the following:\\
$(1)$ \[ \begin{split} \mathfrak{X}_1=& \bigsqcup_
{\begin{smallmatrix}0\leqslant\alpha <p^m,\\ 1\leqslant u \leqslant (p^m-1)/2, p\nmid u\end{smallmatrix}}\widehat{G_l}
\begin{pmatrix} \alpha +u & 0 \\ 0 & \alpha \end{pmatrix}, \
\left|\widehat{G_l}\begin{pmatrix} \alpha +u & 0 \\ 0 & \alpha \end{pmatrix}\right| =p^{2n-2}(p-1)^2, \\
\textit{and}\ \ & \widehat{G_l}\begin{pmatrix} \alpha +u & 0 \\ 0 & \alpha \end{pmatrix}
=\left\{(\mu_{\alpha}\circ\mathrm{det}) \chi\ \Big| \ \chi \in \widehat{G_l}\begin{pmatrix} u & 0 \\0 &0 \end{pmatrix}\right\}.
\end{split}
\]
}\\
\textit{
$(2)$ \[ \begin{split} \mathfrak{X}_2=& \bigsqcup_{0\leqslant\alpha, \epsilon <p^m, 
\left( \frac{\epsilon}{p}\right)=-1}\widehat{G_l}
\begin{pmatrix} \alpha &\epsilon \\ 1 &\alpha \end{pmatrix}, \
\left|\widehat{G_l}\begin{pmatrix} \alpha & \epsilon \\ 1 & \alpha \end{pmatrix}\right| =p^{2n-2}(p^2-1), \\
\textit{and}\ \ & \widehat{G_l}\begin{pmatrix} \alpha & \epsilon \\ 1 & \alpha \end{pmatrix}
=\left\{(\mu_{\alpha}\circ\mathrm{det})\chi\ \Big| \ \chi \in \widehat{G_l}\begin{pmatrix} 0 & \epsilon \\1 &0 \end{pmatrix}\right\},
\end{split}
\]
where $\displaystyle \left(\frac{\epsilon}{p}\right)$ is the Legendre symbol.
}\\
\textit{
$(3)$ \[ \begin{split} \mathfrak{X}_3=& \bigsqcup_{0\leqslant\alpha <p^m, 0\leqslant\beta <p^{m-1}}
\widehat{G_l}
\begin{pmatrix} \alpha &p\beta \\ 1 &\alpha \end{pmatrix}, \
\left|\widehat{G_l}\begin{pmatrix} \alpha & p\beta \\ 1 & \alpha \end{pmatrix}\right| =p^{2n-1}(p-1), \\
\textit{and}\ \ & \widehat{G_l}\begin{pmatrix} \alpha & p\beta \\ 1 & \alpha \end{pmatrix}
=\left\{(\mu_{\alpha}\circ\mathrm{det})\chi\ \Big|\ \chi \in \widehat{G_l}\begin{pmatrix} 0 & p\beta \\1 &0 \end{pmatrix}\right\}.
\end{split}
\]
}\\
\textit{
$(4)$ \begin{gather*}
\mathfrak{X}=\mathfrak{X}_1\bigsqcup\mathfrak{X}_2\bigsqcup\mathfrak{X}_3\ \ \textit{and} \\
|\mathfrak{X}_1|=\frac{1}{2}p^{2l-3}(p-1)^3,
|\mathfrak{X}_2|=\frac{1}{2}p^{2l-3}(p-1)(p^2-1),
|\mathfrak{X}_3|=p^{2l-2}(p-1). 
\end{gather*}
}

\sect{3. Gauss Sums on $GL_2(\Z / p^l\Z)$}

\subsect{3.1.}
Throughout this section, we fix a nontrivial additive character $e$ of $\Z_{p^l}$ and put $e(1)=\zeta_{p^l}^r$. 
Let $\chi$ be a character of $G_l$. Then we define the Gauss sum on $G_l$ associated with $\chi$ by  
\[
\tau_l(\chi) = \sum_{X \in G_l} 
\chi(X) e(\mathrm{Tr}X)=|G_l|\langle \chi, \overline{e\circ\mathrm{Tr}}\rangle_{G_l}.   
\]

Let $H$ be a subgroup of $G_l$ and $\psi$ be a character of $H$. Then, by the Frobenius reciprocity, we have 
\[
\frac{1}{|G_l|}\tau_{l}(\mathrm{ind}_{H}^{G_l}\psi)=\langle \psi, \overline{e\circ\mathrm{Tr}} \rangle_{H}. 
\]

The purpose of this section is to determine explicitly $\tau_l(\chi)$ associated with any irreducible character $\chi$ of $G_l$. 
In the following four subsections, namely $3.2$, $3.3$, $3.4$, and $3.5$, we will evaluate $\tau_l(\chi)$ for every $\chi$ in $\mathfrak{X}_1$, $\mathfrak{X}_2$, 
$\mathfrak{X}_3$, and $\mathfrak{X}_4$ respectively. \\ 

\subsect{3.2. $\mathfrak{X}_1$}

In this subsection, we determine the Gauss sum $\tau_l(\chi)$ for every irreducible character $\chi$ in $\mathfrak{X}_1$.  
Let $0\leqslant\alpha <p^m$, $1\leqslant u \leqslant (p^m-1)/2$ with $p\nmid u$, 
$A=\bigl( \begin{smallmatrix} \alpha +u & 0 \\ 0 & \alpha \end{smallmatrix} \bigr)$, 
$A_0=\bigl( \begin{smallmatrix} u & 0 \\ 0 & 0 \end{smallmatrix} \bigr)$, and 
$T=\mathrm{Stab}_{G_l}\phi_{A_0}$.
Then we have 
\begin{gather*}
T=K_mS=\left\{\begin{pmatrix} a & p^mb \\ p^mc & d  \end{pmatrix}\ 
\Big| \  a,b,c,d \in \Z_{p^l} \right\}\cap G_l ,\\
S=\left\{\begin{pmatrix} a & 0 \\ 0 & d  \end{pmatrix}\ 
\Big| \  a,d \in \Z_{p^l}^{\times} \right\}
\end{gather*}
(cf. \cite{lcw} 3.2.1, 3.3.1 and the proof of Theorem 3.1).
Let
\begin{gather*}
N=\left\{\begin{pmatrix} 1+p^na & p^nb \\ p^mc & 1+p^nd  \end{pmatrix}\ 
\Big| \  a,b,c,d \in \Z_{p^l} \right\}, \\
T_0=\left\{\begin{pmatrix} a & p^nb \\ p^mc & d  \end{pmatrix}\ 
\Big| \  a,b,c,d \in \Z_{p^l} \right\}\cap G_l .
\end{gather*}
Then $N$ is normal in $T$, and $T_0=NS$ (loc. cit.).  
Let $\lambda^{\prime}$ be a multiplicative character of $\Z_{p^l}^{\times}$ such that $\lambda^{\prime}(1+p^n)=\lambda(p^nu)$.
Then we can extend $\phi_{A_0}$ to a character $\phi^{\prime}$ of $N$ by
\[
\phi^{\prime}\begin{pmatrix} 1+p^na & p^nb \\ p^mc & 1+p^nd  \end{pmatrix}=\lambda^{\prime}(1+p^na). 
\]
We have $\mathrm{Stab}_{G_l}\phi^{\prime}=T_0$ (loc. cit.) and we can also extend $\phi^{\prime}$ to characters $\psi_{ij}$ $(0\leqslant i,j <p^{n-1}(p-1))$ of $T_0$ by 
\[
\psi_{ij}\begin{pmatrix} a & p^nb \\ p^mc & d  \end{pmatrix}
=\lambda^{\prime}(a)^{1+p^mi}\lambda^{\prime}(d)^{p^mj}.
\]
By Clifford theory and Theorem $2.1$ $(1)$, the irreducible characters of $G_l$ in $\widehat{G_l}(A)$ are 
$\chi_{ij}=(\mu_{\alpha}\circ\mathrm{det})\mathrm{ind}_{T_0}^{G}\psi_{ij}\ (0\leqslant i,j <p^{n-1}(p-1))$.

Recall that the Gauss sum over $\Z_{p^l}$ associated with a multiplicative character $\nu$ of $\Z_{p^l}^{\times}$ is given by 
\[
g_l(\nu)=\sum_{x\in \Z_{p^l}^{\times}}\nu(x)e(x) . 
\]
The problem of evaluating $g_l(\nu)$ can be reduced to the case where $\nu$ is primitive and $r=1$, and this value was determined explicitly by Odoni as follows (cf. \cite{od} or \cite{bew} Theorem 1.6.2):\\

\subsect{Theorem 3.1.}
\textit{Let $\nu$ be a primitive character of $\Z_{p^l}^{\times}$ normalized such that $\nu(1+p)=\zeta_{p^{l-1}}^{-1}$. If $r=1$, then
\[
g_l(\nu)=
\begin{cases}
p^{l/2}\zeta_{p^l}, & \textrm{if $l=2$}, \\
p^{l/2}\zeta_{p^l}\sqrt{-1}^{(1-p)/2}\zeta_{p}^{(p^2-1)/8}, & \textrm{if $l=3$}, \\
p^{l/2}\zeta_{p^l}^{\sigma}, & \textrm{if $l > 3$ and $l$ is even}, \\
p^{l/2}\zeta_{p^l}^{\sigma}\sqrt{-1}^{(1-p)/2}, & \textrm{if $l > 3$ and $l$ is odd},
\end{cases}
\]
where $\sigma$ is a $p$-adic integer defined by
\[
\sigma= \frac{p}{\log (1+p)}\left( 1-\log \left(\frac{p}{\log (1+p)} \right)\right).
\]
}\\

\subsect{Theorem 3.2.}
\textit{With the notation as above, we have
\[
\frac{1}{\chi_{ij}(I)}\tau_l(\chi_{ij})
=p^l g_l\left(\mu_\alpha(\lambda^{\prime})^{1+p^mi})g_l(\mu_\alpha(\lambda^{\prime})^{p^mj}\right).
\]
}\\
\noindent
\textit{Proof.}
\[
\begin{split}
& \frac{1}{\chi_{ij}(I)}\tau_l(\chi_{ij})=|T_0|\langle (\mu_{\alpha}\circ\mathrm{det})\psi_{ij}, \overline{e\circ\mathrm{Tr}} \rangle_{T_0}\\
=& \sum_{\begin{smallmatrix}0\leqslant a,d <p^l,\\ p\nmid ad\end{smallmatrix}}
\sum_{\begin{smallmatrix}0\leqslant b <p^m,\\ 0\leqslant c <p^n\end{smallmatrix}}\mu_{\alpha}(ad)
\lambda^{\prime}(a)^{1+p^mi}\lambda^{\prime}(d)^{p^mj}e(a+d)\\
=& p^l\left(\sum_{0\leqslant a <p^l, p\nmid a}\mu_{\alpha}(a)\lambda^{\prime}(a)^{1+p^mi}e(a)\right)
\left(\sum_{0\leqslant d <p^l, p\nmid d}\mu_{\alpha}(d)\lambda^{\prime}(d)^{p^mj}e(d)\right)\\
=& p^l g_l\left(\mu_\alpha(\lambda^{\prime})^{1+p^mi})g_l(\mu_\alpha(\lambda^{\prime})^{p^mj}\right),
\end{split}
\]
which proves the assertion of the theorem.
\hfill
$\Box$\\

\subsect{3.3. $\mathfrak{X}_2$}

In this subsection, we determine the Gauss sum $\tau_l(\chi)$ for every irreducible character $\chi$ in $\mathfrak{X}_2$.  
Let $0\leqslant\alpha, \epsilon <p^m$ with $\left( \frac{\epsilon}{p}\right)=-1$, 
$A=\bigl( \begin{smallmatrix} \alpha & \epsilon \\ 1 & \alpha \end{smallmatrix} \bigr)$, 
$A_0=\bigl( \begin{smallmatrix} 0 & \epsilon \\ 1 & 0 \end{smallmatrix} \bigr)$, 
and $T=\mathrm{Stab}_{G_l}\phi_{A_0}$.
Then we have
\begin{gather*}
T=K_mS=\left\{\begin{pmatrix} a & \epsilon b+p^mc \\ b & a+p^md  \end{pmatrix}\ 
\Big| \  a,b,c,d \in \Z_{p^l} \right\}\cap G_l ,\\
S=\left\{\begin{pmatrix} a & \epsilon b \\ b & a  \end{pmatrix}\ 
\Big| \  a,b \in \Z_{p^l}  \right\}\cap G_l
\end{gather*}
(cf. \cite{lcw} 3.2.2 and 3.3.2).\\
Let $s_1= \bigl(\begin{smallmatrix}1+p & \\ & 1+p \end{smallmatrix}\bigr)$, 
$s_2= \bigl(\begin{smallmatrix}1 & p\epsilon \\ p & 1\end{smallmatrix}\bigr)$. We have
\[
K_1 \cap S=\langle s_1 \rangle \times \langle s_2 \rangle ,\ K_n \cap S= \langle s_1^{p^{n-1}} \rangle \times \langle s_2^{p^{n-1}} \rangle .
\] 
Moreover, there exists $s_3 \in S$ such that
\[
S=(K_1 \cap S)\times \langle s_3  \rangle,\ s_3^{p+1} \in Z(G_l),
\]
where $Z(G_l)$ denotes the center of $G_l$ (cf. \cite{lcw} Lemma 5.1).
Hence the coset representatives of $K_n \cap S$ in $S$ are given by $s^k=s_1^{k_1}s_2^{k_2}s_3^{k_3} \ (k \in \Omega)$, 
where $\Omega=\{k=(k_1, k_2, k_3) \ | \ 0\leqslant k_1, k_2 < p^{n-1}, 0 \leqslant k_3 < p^2-1 \}$. 
Also we define $\delta \in\Z_{p^l}^{\times}$ by $p^n \delta = 2\sum_{1\leqslant t \leqslant p^{n-1}, t:\mathrm{odd}}
\bigl(\begin{smallmatrix} p^{n-1} \\ t\end{smallmatrix}\bigr)p^t\epsilon^{(t+1)/2}$.

We evaluate the Gauss sum $\tau_l(\chi)$ for $\chi$ in $\mathfrak{X}_2$ in two cases depending on whether $l$ is even or odd. \\ 

\subsect{3.3.1. Even case, i.e. $l=2m$}

The number of extensions of $\phi_{A_0}$ to $T$ is $|T:K_m|=|\Omega|$, and
\begin{gather*}
\phi_{A_0}(s_1^{p^{m-1}})=1,\ \phi_{A_0}(s_3^{p^2-1})=1, \\
\phi_{A_0}(s_2^{p^{m-1}})=\lambda(p^m\delta)=\zeta_{p^m}^{\delta}.
\end{gather*}
Hence for every $i\in \Omega$, there exists a unique extension $\psi_ i$ of $\phi_{A_0}$ such that
\[
\psi_i(s_1)=\zeta_{p^{m-1}}^{i_1},\ \psi_i(s_2)=\zeta_{p^{l-1}}^{\delta +p^m i_2},\ \psi_i(s_3)=\zeta_{p^2-1}^{i_3}.
\]
By Clifford theory and Theorem $2.1$ $(2)$,  
the irreducible characters of $G_l$ in $\widehat{G_l}(A)$ are $\chi_i=(\mu_\alpha\circ\det)\mathrm{ind}_{T}^{G_l}\psi_i \ (i \in \Omega)$.\\

\subsect{Theorem 3.3.}
\textit{With the notation as above, we have the following.\\ 
$(1)$ If $p | r$, then $\tau_l(\chi_i)=0$.\\
$(2)$ If $p \nmid r$, then
\[
\frac{1}{\chi_i(I)}\tau_l(\chi_i)=p^{2l}\mu_{\alpha}(\det s^h)e(\mathrm{Tr}s^h)\zeta_{p^{m-1}}^{i_1h_1}\zeta_{p^{l-1}}^{(\delta +p^m i_2)h_2}\zeta_{p^2-1}^{i_3h_3},
\]
where $h$ is a unique element of $\Omega$ such that
$\bigl(\begin{smallmatrix}-\alpha r^{-1}& -\epsilon r^{-1} \\ -r^{-1}& -\alpha r^{-1}\end{smallmatrix}\bigr)\in (K_m \cap S)s^h$.
}\\
\noindent
\textit{Proof.}
For $k \in \Omega$, let $s^k=\bigl(\begin{smallmatrix}x_k & \epsilon y_k \\ y_k & x_k\end{smallmatrix}\bigl)$. Then
\[
\begin{split}
& \frac{1}{\chi_{i}(I)}\tau_l(\chi_{i})=|T|\langle (\mu_{\alpha}\circ\mathrm{det})\psi_{i}, \overline{e\circ\mathrm{Tr}} \rangle_{T}\\
=& \sum_{0\leqslant a,b,c,d <p^m}\sum_{k \in \Omega}\mu_{\alpha}((1+p^ma)(1+p^md))\mu_{\alpha}(\det s^k)\lambda (p^m\epsilon c +p^mb)\psi_i(s^k)\\
& \cdot e(\mathrm{Tr}s^k)
e( p^m\mathrm{Tr}\bigl( \begin{smallmatrix}a & b \\ c& d \end{smallmatrix}\bigr)
\bigl(\begin{smallmatrix}x_k & \epsilon y_k \\ y_k & x_k\end{smallmatrix}\bigl) )\\
=& \sum_{k \in \Omega}\mu_{\alpha}(\det s^k)e(\mathrm{Tr}s^k)\psi_i(s^k) \\
& \cdot \left( \sum_{0\leqslant a <p^m}\lambda (p^m \alpha a) e(p^m x_k a)\right)^2 \left(\sum_{0\leqslant b <p^m}\lambda (p^m b)e(p^m y_k b)\right)^2 \\
=& p^{2l} \sum_{\begin{smallmatrix}k \in \Omega , \\ \alpha +x_kr \equiv 0\ (\mathrm{mod}\ p^m), \\1+y_kr \equiv 0\ (\mathrm{mod}\ p^m) \end{smallmatrix}}
\mu_{\alpha}(\det s^k)e(\mathrm{Tr}s^k)\psi_i(s^k) .
\end{split}
\]
\noindent
Thus, $\tau_l(\chi)=0$ if $p|r$.

Assume $p\nmid r$. Then there exists a unique element $h$ of $\Omega$ such that
$\bigl(\begin{smallmatrix}-\alpha r^{-1}& -\epsilon r^{-1} \\ -r^{-1}& -\alpha r^{-1}\end{smallmatrix}\bigr)\in (K_m \cap S)s^h$.
For every $k\in\Omega$, 
$\alpha +x_kr \equiv 0\ (\mathrm{mod}\ p^m)$ and $1+y_kr \equiv 0\ (\mathrm{mod}\ p^m)$ if and only if 
$\bigl(\begin{smallmatrix}-\alpha r^{-1}& -\epsilon r^{-1} \\ -r^{-1}& -\alpha r^{-1}\end{smallmatrix}\bigr)\in (K_m \cap S)s^k$, namely $k=h$. 
Hence 
\[
\begin{split}
& \frac{1}{\chi_{i}(I)}\tau_l(\chi_{i})= p^{2l}\mu_{\alpha}(\det s^h)e(\mathrm{Tr}s^h)\psi_i(s^h)\\
=& p^{2l}\mu_{\alpha}(\det s^h)e(\mathrm{Tr}s^h)\zeta_{p^{m-1}}^{i_1h_1}\zeta_{p^{l-1}}^{(\delta +p^m i_2)h_2}\zeta_{p^2-1}^{i_3h_3},
\end{split}
\]
which proves the assertion of the theorem.
\hfill
$\Box$\\

\subsect{3.3.2. Odd case, i.e. $l=2m+1$}

Let $L=K_{m+1}S$. The number of extensions of $\phi_{A_0}$ to $L$ is $|T:K_{m+1}|=|\Omega|$, and
\begin{gather*}
\phi_{A_0}(s_1^{p^{m}})=1,\ \phi_{A_0}(s_3^{p^2-1})=1, \\
\phi_{A_0}(s_2^{p^{m}})=\lambda(p^{m+1}\delta)=\zeta_{p^m}^{\delta}.
\end{gather*}
Hence for every $i\in \Omega$, there exists a unique extension $\phi_ i$ of $\phi_{A_0}$ such that
\[
\phi_i(s_1)=\zeta_{p^{m}}^{i_1},\ \phi_i(s_2)=\zeta_{p^{l-1}}^{\delta +p^m i_2},\ \phi_i(s_3)=\zeta_{p^2-1}^{i_3}.
\]
Let
\[
N_j=K_jZ(G_l)(K_1\cap S)=\left\{\begin{pmatrix} a & p\epsilon b+p^jc \\ pb & a+p^jd  \end{pmatrix}\ 
\Big| \  a,b,c,d \in \Z_{p^l} \right\}\cap G_l,\ j=m,m+1.
\]
Then $N_j$ is normal in $T$ (cf. \cite{lcw} 3.3.2). Let $\phi_i^{\prime} = \phi_i |_{N_{m+1}}$.  
Then $\phi_i^{\prime}$ is stable under $T$. 
Let $H=N_{m+1}\langle \bigl(\begin{smallmatrix}1+p^m &  \\ & 1\end{smallmatrix}\bigr) \rangle$. 
Since $N_m/N_{m+1}$ is abelian, $H$ is normal in $N_m$. We can extend $\phi_i^{\prime}$ to a linear character $\phi_i^{\prime \prime}$ of $H$ 
such that $\phi_i^{\prime \prime}$ is trivial on $\langle \bigl(\begin{smallmatrix}1+p^m &  \\ & 1\end{smallmatrix}\bigr) \rangle$. 
It is easy to find an element in $N_m$ that does not stabilize $\phi_i^{\prime \prime}$ and since the index of $H$ in $N_m$ is $p$, 
we have $\mathrm{Stab}_{N_m} \phi_i^{\prime \prime}=H$. Then Clifford theory implies that $\mathrm{ind}_{H}^{N_m}\phi_i^{\prime \prime}$ is irreducible. 
Since $\phi_i^{\prime}$ is stable under $T$, we have $(\mathrm{ind}_{H}^{N_m}\phi_i^{\prime \prime}) |_{N_{m+1}}=p\phi_i^{\prime}$. 
Hence $(1/p)\mathrm{ind}_{N_{m+1}}^{N_m}\phi_i^{\prime}= \mathrm{ind}_{H}^{N_m}\phi_i^{\prime \prime} \in \widehat{N_m}$.\\
 
\subsect{Lemma 3.4.}
\textit{Let $\psi_i=(1/p)\mathrm{ind}_{N_{m+1}}^{T}\phi^{\prime}_i-\mathrm{ind}_{L}^{T}\phi_i$ for $i\in\Omega$. Then 
$\psi_i$, for $i\in \Omega$, are distinct irreducible characters of $T$ and satisfy $\psi_i|_{K_{m+1}}=p\phi_{A_0}$.
}
\\
\noindent
\textit{Proof.}
Since $\phi_i^{\prime}$ is stable under $T$, we have
\[
(\mathrm{ind}_{N_{m+1}}^{T}\phi^{\prime}_i)|_{N_{m+1}}=p^2(p+1)\phi^{\prime}_i,\ (\mathrm{ind}_{L}^{T}\phi_i )|_{N_{m+1}} = p^2\phi^{\prime}_i.
\]
Hence we have $\psi_i|_{K_{m+1}}=p\phi_{A_0}$ and by the Frobenius reciprocity
\[
\langle \mathrm{ind}_{N_{m+1}}^{T}\phi_i^{\prime}, \mathrm{ind}_{N_{m+1}}^{T}\phi_i^{\prime}\rangle_{T}=p^2(p+1),\ 
\langle \mathrm{ind}_{N_{m+1}}^{T}\phi_i^{\prime}, \mathrm{ind}_{L}^{T}\phi_i\rangle_{T}=p^2.
\]
We can pick the right coset representatives of $L$ in $T$ to be
$E_{cd}=\bigl( \begin{smallmatrix} 1 & p^mc \\ 0 & 1+p^md \end{smallmatrix} \bigr)$ $(0\leqslant c,d <p)$.
If $(c,d)\not=(0,0)$, then $L\cap E_{cd}^{-1}LE_{cd}=N_{m+1}$, so the number of right cosets in $LE_{cd}L$ is $|L:N_{m+1}|=p+1$.
Hence there exists a suitable subset $\Gamma$ of $\{(c,d) \ |\ 0\leqslant c,d <p \}$ with $(0,0) \in \Gamma$ and $|\Gamma|= p$ such that 
$E_{cd}$ $((c,d)\in \Gamma)$ are the double coset representatives of $L$ in $T$. 
By the Mackey theorem, we have
\[
\begin{split}
\langle \mathrm{ind}_{L}^{T}\phi_i, \mathrm{ind}_{L}^{T}\phi_i\rangle_{T}=& \sum_{(c,d)\in \Gamma}\langle
\phi_i^{E_{cd}}|_{L\cap E_{cd}^{-1}LE_{cd}}, \phi_i |_{L\cap E_{cd}^{-1}LE_{cd}}\rangle_{L\cap E_{cd}^{-1}LE_{cd}}\\
=& 1+ \sum_{(c,d)\in \Gamma, (c,d)\not=(0,0)}\langle
\phi_i^{E_{cd}}|_{N_{m+1}}, \phi_i^{\prime}\rangle_{N_{m+1}}\\
=& 1+ \sum_{(c,d)\in \Gamma, (c,d)\not=(0,0)}\langle \phi_i^{\prime}, \phi_i^{\prime}\rangle_{N_{m+1}}=p. 
\end{split}
\]
Therefore 
\[
\langle \psi_i,\psi_i \rangle_{T}= \frac{1}{p^2}\langle \mathrm{ind}_{N_{m+1}}^{T}\phi_i^{\prime}, \mathrm{ind}_{N_{m+1}}^{T}\phi_i^{\prime}\rangle_{T}
-\frac{2}{p}\langle \mathrm{ind}_{N_{m+1}}^{T}\phi_i^{\prime}, \mathrm{ind}_{L}^{T}\phi_i\rangle_{T}+\langle \mathrm{ind}_{L}^{T}\phi_i, 
\mathrm{ind}_{L}^{T}\phi_i\rangle_{T}=1.
\]
Since $\psi_i(I)=p$, we conclude that $\psi_i$ is irreducible. 

Suppose $\psi_i=\psi_j$ for $i,j\in\Omega$.
Restricting to $N_{m+1}$ on both sides, we have $\phi_i^{\prime}=\phi_j^{\prime}$. 
Hence
\[
\begin{split}
1=& \langle \psi_i,\psi_j \rangle_{T}\\
=& \frac{1}{p^2}\langle \mathrm{ind}_{N_{m+1}}^{T}\phi_i^{\prime}, \mathrm{ind}_{N_{m+1}}^{T}\phi_i^{\prime}\rangle_{T}
-\frac{1}{p}\langle \mathrm{ind}_{N_{m+1}}^{T}\phi_i^{\prime}, \mathrm{ind}_{L}^{T}\phi_i\rangle_{T}\\
& -\frac{1}{p}\langle \mathrm{ind}_{N_{m+1}}^{T}\phi_j^{\prime}, \mathrm{ind}_{L}^{T}\phi_j\rangle_{T}
+\langle \mathrm{ind}_{L}^{T}\phi_i, \mathrm{ind}_{L}^{T}\phi_j\rangle_{T}\\
=& -p+1+\langle \phi_i, \phi_j\rangle_{L}+\sum_{(c,d)\in \Gamma, (c,d)\not=(0,0)}\langle \phi_i^{E_{cd}}|_{N_{m+1}}, \phi_i^{\prime}\rangle_{N_{m+1}}\\
=& \langle \phi_i, \phi_j\rangle_{L}, 
\end{split}
\]
which yields $i=j$. The proof is complete. 
\hfill
$\Box$\\

By this Lemma, Clifford theory, and Theorem $2.1$ $(2)$, 
the irreducible characters of $G_l$ in $\widehat{G_l}(A)$ are $\chi_i= (\mu_{\alpha}\circ \det)\mathrm{ind}_{T}^{G_l}\psi_i\ (i\in \Omega) $.\\

\subsect{Theorem 3.5.}
\textit{With the notation as above, we have the following.\\
$(1)$ If $p | r$, then $\tau_l(\chi_i)=0$.\\
$(2)$ If $p \nmid r$, then
\[
\frac{1}{\chi_i(I)}\tau_l(\chi_i)
=-p^{2l-1}\sum_{k}\mu_{\alpha}(\det s^k)e(\mathrm{Tr}s^k)\zeta_{p^{m}}^{i_1k_1}\zeta_{p^{l-1}}^{(\delta +p^m i_2)k_2}\zeta_{p^2-1}^{i_3h_3},
\]
where $h$ is a unique element of $\Omega$ which satisfies $0\leqslant h_1,h_2 <p^{m-1}$ and
$\bigl(\begin{smallmatrix}-\alpha r^{-1}& -\epsilon r^{-1} \\ -r^{-1}& -\alpha r^{-1}\end{smallmatrix}\bigr)\in (K_m \cap S)s^h$, and
the summation is over all elements $k\in\Omega$ such that $k_1 \equiv h_1$ $(\mathrm{mod}\ p^{m-1})$, $k_2 \equiv h_2\ (\mathrm{mod}\ p^{m-1})$, $k_3=h_3$.\\
}
\noindent
\textit{Proof.} By the Frobenius reciprocity, we have
\[
\frac{1}{|G_l|}\tau_l(\chi)=\frac{1}{p}\langle (\mu_{\alpha}\circ\mathrm{det})\phi_{i}^{\prime}, \overline{e\circ\mathrm{Tr}} \rangle_{N_{m+1}}
-\langle (\mu_{\alpha}\circ\mathrm{det})\phi_{i}, \overline{e\circ\mathrm{Tr}} \rangle_{L}.
\]
Let us show $\langle (\mu_{\alpha}\circ\mathrm{det})\phi_{i}^{\prime}, \overline{e\circ\mathrm{Tr}} \rangle_{N_{m+1}}=0$.
Let $s^k=\bigl(\begin{smallmatrix}x_k & \epsilon y_k \\ y_k & x_k\end{smallmatrix}\bigl)$ for $k \in \Omega$ and 
let $\Omega_0$ be the set of all elements $k$ of $\Omega$ with $p+1|k_3$. Since
\[
Z(G_l)(K_1\cap S)=\langle s_1 \rangle \times\langle s_2 \rangle \times\langle s_3^{p+1} \rangle, 
K_{m+1}\cap Z(G_l)(K_1\cap S)=\langle s_1^{p^m} \rangle \times\langle s_2^{p^m} \rangle ,  
\]
$s^k$ $(k\in\Omega_{0})$ are the coset representatives of $K_{m+1}$ in $N_{m+1}$.
Hence
\[
\begin{split}
&|N_{m+1}|\langle (\mu_{\alpha}\circ\mathrm{det})\phi_{i}^{\prime}, \overline{e\circ\mathrm{Tr}} \rangle_{N_{m+1}}\\
=& \sum_{0\leqslant a,b,c,d <p^m}\sum_{k \in \Omega_0}\mu_{\alpha}((1+p^{m+1}a)(1+p^{m+1}d))\mu_{\alpha}(\det s^k)\lambda (p^{m+1}\epsilon c +p^{m+1}b)\phi^{\prime}_i(s^k)\\
& \cdot e(\mathrm{Tr}s^k)
e(p^{m+1}\mathrm{Tr}\bigl( \begin{smallmatrix}a & b \\ c& d \end{smallmatrix}\bigr)
\bigl(\begin{smallmatrix}x_k & \epsilon y_k \\ y_k & x_k\end{smallmatrix}\bigl) )\\
=& \sum_{k \in \Omega_0}\mu_{\alpha}(\det s^k)e(\mathrm{Tr}s^k)\phi^{\prime}_i(s^k) \\
& \cdot \left( \sum_{0\leqslant a <p^{m}}\lambda (p^{m+1} \alpha a) e(p^{m+1} x_k a)\right)^2 \left(\sum_{0\leqslant b <p^{m}}\lambda (p^{m+1} b)e(p^{m+1} y_k b)\right)^2 .
\end{split}
\] 
Since $p|y_k$ for all $k\in\Omega_{0}$, we have
\[
\sum_{0\leqslant b <p^{m}}\lambda (p^{m+1} b)e(p^{m+1} y_k b)=0, 
\]
which yields $\langle (\mu_{\alpha}\circ\mathrm{det})\phi_{i}^{\prime}, \overline{e\circ\mathrm{Tr}} \rangle_{N_{m+1}}=0$.

Similarly we have
\[
\begin{split}
&|L|\langle (\mu_{\alpha}\circ\mathrm{det})\phi_{i}, \overline{e\circ\mathrm{Tr}} \rangle_{L}\\
=& \sum_{k \in \Omega}\mu_{\alpha}(\det s^k)e(\mathrm{Tr}s^k)\phi_i(s^k) \\
& \cdot \left( \sum_{0\leqslant a <p^{m}}\lambda (p^{m+1} \alpha a) e(p^{m+1} x_k a)\right)^2 \left(\sum_{0\leqslant b <p^{m}}\lambda (p^{m+1} b)e(p^{m+1} y_k b)\right)^2 \\
=& p^{2l-2} \sum_{\begin{smallmatrix}k \in \Omega , \\ \alpha +x_kr \equiv 0\ (\mathrm{mod}\ p^m), \\1+y_kr \equiv 0\ (\mathrm{mod}\ p^m) \end{smallmatrix}}
\mu_{\alpha}(\det s^k)e(\mathrm{Tr}s^k)\phi_i(s^k) .
\end{split}
\] 
\noindent
Thus, $\tau_l(\chi)=0$ if $p|r$.

Assume $p\nmid r$. Let $\Omega_1$ be the set of all elements $k$ of $\Omega$ with $0\leqslant k_1, k_2 <p^{m-1}$.
Since $K_m\cap S=\langle s_1^{p^{m-1}} \rangle \times\langle s_2^{p^{m-1}}\rangle$, 
$s^k\ (k\in\Omega_1)$ are the coset representatives of $K_m\cap S$ in $S$.
Hence, there exists a unique element $h$ of $\Omega_1$ such that
$\bigl(\begin{smallmatrix}-\alpha r^{-1}& -\epsilon r^{-1} \\ -r^{-1}& -\alpha r^{-1}\end{smallmatrix}\bigr)\in (K_m \cap S)s^h$. 
For $k\in\Omega$, $\alpha +x_kr \equiv 0\ (\mathrm{mod}\ p^m)$ and $1+y_kr \equiv 0\ (\mathrm{mod}\ p^m)$ if and only if 
$\bigl(\begin{smallmatrix}-\alpha r^{-1}& -\epsilon r^{-1} \\ -r^{-1}& -\alpha r^{-1}\end{smallmatrix}\bigr)\in (K_m \cap S)s^k$, namely
$k_1 \equiv h_1\ (\mathrm{mod}\ p^{m-1}), k_2 \equiv h_2\ (\mathrm{mod}\ p^{m-1}), k_3=h_3$.
Therefore
\[
\begin{split}
& |L|\langle (\mu_{\alpha}\circ\mathrm{det})\phi_{i}, \overline{e\circ\mathrm{Tr}} \rangle_{L}\\
=& p^{2l-2}\sum_{\begin{smallmatrix}k\in\Omega,\\ k_1 \equiv h_1\ (\mathrm{mod}\ p^{m-1}),\\ k_2 \equiv h_2\ (\mathrm{mod}\ p^{m-1}),\\ k_3=h_3\end{smallmatrix}}
\mu_{\alpha}(\det s^k)e(\mathrm{Tr}s^k)\zeta_{p^{m}}^{i_1k_1}\zeta_{p^{l-1}}^{(\delta +p^m i_2)k_2}\zeta_{p^2-1}^{i_3h_3},
\end{split}
\]
which proves the assertion of the theorem.  
\hfill
$\Box$\\

\subsect{3.4. $\mathfrak{X}_3$}

In this subsection, we determine the Gauss sum $\tau_l(\chi)$ for every irreducible character $\chi$ in $\mathfrak{X}_3$.  
Let $0\leqslant \alpha <p^{m},0\leqslant \beta <p^{m-1}$, 
$A=\bigl( \begin{smallmatrix} \alpha & p\beta \\ 1 & \alpha \end{smallmatrix} \bigr)$, 
$A_0=\bigl( \begin{smallmatrix} 0 & p\beta \\ 1 & 0 \end{smallmatrix} \bigr)$, 
and $T=\mathrm{Stab}_{G_l}\phi_{A_0}$.
Then we have 
\begin{gather*}
T=K_mS=\left\{\begin{pmatrix} a & p\beta b+p^mc \\ b & a+p^md  \end{pmatrix} 
\Big| \  a,b,c,d \in \Z_{p^l} \right\}\cap G_l, \\
S=\left\{\begin{pmatrix} a & p\beta b \\ b & a  \end{pmatrix}\ 
\Big| \  a,b \in \Z_{p^l}  \right\}\cap G_l
\end{gather*}
(cf. \cite{lcw} 3.2.3 and 3.3.3). Let
\[
N=\left\{\begin{pmatrix} 1+p^ma & p^nb \\ p^mc & 1+p^md  \end{pmatrix}\ 
\Big| \  a,b,c,d \in \Z_{p^l} \right\}.
\]
Then $N$ is normal in $T$ 
and we can extend $\phi_{A_0}$ to characters $\phi_i$ $(i=(i_1, i_2),\ 0\leqslant i_1,i_2 <p^{n-m})$ of $N$ defined by
\[
\phi_{i}\begin{pmatrix} 1+p^ma & p^nb \\ p^mc & 1+p^md  \end{pmatrix}
=\lambda(p^{2m}(i_1a+i_2d))\lambda(p^{m+1}\beta c+p^nb).
\]
Let $T_0=\mathrm{Stab}_{T}\phi_{i}$. Then we have 
\[
T_0=NS=\left\{\begin{pmatrix} a & p\beta b+p^nc \\ b & a+p^md  \end{pmatrix}\ 
\Big| \  a,b,c,d \in \Z_{p^l} \right\}\cap G_l.
\]
\\
Let $\gamma$ be a generator of $\Z_{p^l}^{\times}$ such that $\gamma^{p^{m-1}(p-1)}=1+p^m$, and let
$s_1= \bigl(\begin{smallmatrix}\gamma & \\ & \gamma \end{smallmatrix}\bigr), 
s_2= \bigl(\begin{smallmatrix}1 & p^2\beta \\ p & 1\end{smallmatrix}\bigr), 
s_3= \bigl(\begin{smallmatrix}1 & p\beta \\ 1 & 1\end{smallmatrix}\bigr).$
Then we have
\begin{gather*}
Z(G_l)(K_1 \cap S)=\langle s_1 \rangle \times \langle s_2 \rangle ,\ 
S = (N \cap S)\langle s_1 \rangle \langle s_2 \rangle \langle s_3 \rangle, \\
K_m \cap S= N \cap S= \langle s_1^{p^{m-1}(p-1)} \rangle \times \langle s_2^{p^{m-1}} \rangle,
\end{gather*}
and $s_3^p=s_1^{\sigma_1}s_2^{\sigma_2}$ for some integers $\sigma_1$ and $\sigma_2$.
Hence the coset representatives of $N$ in $T_0$ are given by $s^k=s_1^{k_1}s_2^{k_2}s_3^{k_3}$ $( k\in\Omega )$, 
where $\Omega=\{k=(k_1, k_2, k_3) \ | \ 0\leqslant k_1 < p^{m-1}(p-1), 0 \leqslant k_2 < p^{m-1}, 0\leqslant k_3< p \}$. 
Also we define $\delta \in\Z_{p^l}^{\times}$ by $p^{m+1} \delta = 2\sum_{1\leqslant t \leqslant p^{m-1}, t:\mathrm{odd}}
\bigl(\begin{smallmatrix} p^{m-1} \\ t\end{smallmatrix}\bigr)p^t(p\beta)^{(t+1)/2}$.
Then we have
\begin{gather*}
\phi_{i}(s_1^{p^{m-1}(p-1)})=\phi_{i}((1+p^m)I)=\lambda(p^{2m}(i_1+i_2))= \zeta_{p^{n-m}}^{i_1+i_2},\\
\phi_{i}(s_2^{p^{m-1}})=\lambda(p^{m+1}\delta)=\zeta_{p^{n-1}}^{\delta}.
\end{gather*}
Since the number of extensions of $\phi_{i}$ to $T_0$ is $|T_0:N|=|\Omega|$, for every $j\in\Omega$, 
there exists a unique extension $\psi_{ij}$ of $\phi_{i}$ such that 
\begin{gather*}
\psi_{ij}(s_1)=\zeta_{p^{n-1}(p-1)}^{i_1+i_2+p^{n-m}j_1},\ \psi_{ij}(s_2)=\zeta_{p^{l-2}}^{\delta +p^{n-1} j_2},\\
\psi_{ij}(s_3)=\zeta_{p^{n}(p-1)}^{(i_1+i_2+p^{n-m}j_1)\sigma_1}\zeta_{p^{l-1}}^{(\delta +p^{n-1} j_2)\sigma_2}\zeta_{p}^{j_3}.
\end{gather*}
Therefore, by Clifford theory and Theorem $2.1$ $(3)$, 
the irreducible characters of $G_l$ in $\widehat{G_l}(A)$ are 
$\chi_{ij}=(\mu_\alpha\circ\det)\mathrm{ind}_{T_0}^{G_l}\psi_{ij} \ (0\leqslant i_1, i_2<p^{n-m},\ j \in \Omega)$.\\

\subsect{Theorem 3.6.}
\textit{With the notation as above, we have the following.\\
$(1)$ If $p | r$ or $p | \alpha$, then $\tau_l(\chi_{ij})=0$.\\
$(2)$ If $p \nmid r$ and $p \nmid \alpha$, then
\[
\begin{split}
\frac{1}{\chi_{ij}(I)}\tau_l(\chi_{ij})
=& p^{l+2m}\mu_{\alpha}(\det s^h)e(\mathrm{Tr}s^h)\zeta_{p^{n}(p-1)}^{(i_1+i_2+p^{n-m}j_1)(ph_1+\sigma_1 h_3)}\zeta_{p^{l-1}}^{(\delta +p^{n-1} j_2)(ph_2+\sigma_2 h_3)}\zeta_{p}^{j_3h_3}\\
& \cdot \Big(\sum_{0\leqslant t_1< p^{n-m}}\mu_{\alpha}(1+p^mt_1)\lambda(p^{2m}i_1t_1)e(p^m{x_h}t_1)\Big)\\
& \cdot \Big(\sum_{0\leqslant t_2 < p^{n-m}}\mu_{\alpha}(1+p^mt_2)\lambda(p^{2m}i_2t_2)e(p^m{x_h}t_2)\Big),
\end{split}
\]
where $s^k=\bigl(\begin{smallmatrix}x_k & p\beta y_k \\ y_k & x_k\end{smallmatrix}\bigl)$ for $k\in\Omega$ and
$h$ is a unique element of $\Omega$ which satisfies 
$\bigl(\begin{smallmatrix}-\alpha r^{-1}& -p\beta r^{-1} \\ -r^{-1}& -\alpha r^{-1}\end{smallmatrix}\bigr)\in (K_m \cap S)s^h$. \\
}
\noindent
\textit{Proof.}
We define a subgroup $M$ of $G_l$ by
\[
M=\left\{\begin{pmatrix} 1+p^na & p^nb \\ p^mc & 1+p^nd  \end{pmatrix}\ 
\Big| \  a,b,c,d \in \Z_{p^l} \right\},
\]
and pick the coset representatives of $M$ in $N$ to be
\[
X_t=\bigl(\begin{smallmatrix}1+p^mt_1 &  \\ & 1+p^mt_2 \end{smallmatrix}\bigr),\ t=(t_1,t_2),\ 0\leqslant t_1, t_2 <p^{n-m}.
\]
Then we have
\[
\begin{split}
& \frac{1}{\chi_{ij}(I)}\tau_l(\chi_{ij})=|T_0|\langle (\mu_{\alpha}\circ\mathrm{det})\psi_{ij}, \overline{e\circ\mathrm{Tr}} \rangle_{T_0}\\
=& \sum_{\begin{smallmatrix}0\leqslant a,b,d <p^m,\\ 0\leqslant c <p^n\end{smallmatrix}}\sum_{\begin{smallmatrix}0\leqslant t_1, t_2 <p^{n-m}, \\ k \in \Omega \end{smallmatrix}}
\mu_{\alpha}((1+p^{n}a)(1+p^{n}d))\mu_{\alpha}((1+p^mt_1)(1+p^mt_2))\mu_{\alpha}(\det s^k)\\
& \cdot \lambda (p^{m+1}\beta c +p^{n}b) \lambda (p^{2m}(i_1t_1+i_2t_2))\psi_{ij}(s^k) e(\mathrm{Tr}(I+\bigl( \begin{smallmatrix}p^na & p^nb \\ p^mc& p^nd \end{smallmatrix}\bigr))X_ts^k)\\
=& \sum_{\begin{smallmatrix}0\leqslant t_1, t_2 <p^{n-m}, \\ k \in \Omega \end{smallmatrix}}\mu_{\alpha}((1+p^mt_1)(1+p^mt_2))\mu_{\alpha}(\det s^k)
\lambda (p^{2m}(i_1t_1+i_2t_2))e(\mathrm{Tr}X_ts^k)\psi_{ij}(s^k) \\
& \cdot \left( \sum_{0\leqslant a <p^{m}}\lambda (p^{n} \alpha a) e(p^{n} x_k a)\right)^2 \left(\sum_{0\leqslant b <p^{m}}\lambda (p^{n} b)e(p^{n} y_k b)\right)\\ 
& \cdot \left(\sum_{0\leqslant c <p^{n}}\lambda (p^{m+1}\beta c)e(p^{m+1} y_k \beta c)\right) \\
=& p^{l+2m} \sum_{0\leqslant t_1, t_2 <p^{n-m}} \sum_{\begin{smallmatrix}k \in \Omega , \\ \alpha +x_kr \equiv 0\ (\mathrm{mod}\ p^m), \\1+y_kr \equiv 0\ (\mathrm{mod}\ p^m) \end{smallmatrix}}
\mu_{\alpha}((1+p^mt_1)(1+p^mt_2))\lambda(p^{2m}(i_1t_1+i_2t_2))\\
& \cdot \mu_{\alpha}(\det s^k)e(\mathrm{Tr}s^k) e(p^mx_k(t_1+t_2))\psi_{ij}(s^k).
\end{split}
\]
Thus, $\tau_l(\chi)=0$ if $p|r$ or $p|\alpha$.

Assume $p\nmid r$ and $p\nmid \alpha$.
Then there exists a unique element $h\in\Omega$ such that
$\bigl(\begin{smallmatrix}-\alpha r^{-1}& -p\beta r^{-1} \\ -r^{-1}& -\alpha r^{-1}\end{smallmatrix}\bigr)\in (K_m \cap S)s^h$. 
For every $k\in\Omega$, 
$\alpha +x_kr \equiv 0\ (\mathrm{mod}\ p^m)$ and $1+y_kr \equiv 0\ (\mathrm{mod}\ p^m)$ if and only if 
$\bigl(\begin{smallmatrix}-\alpha r^{-1}& -p\beta r^{-1} \\ -r^{-1}& -\alpha r^{-1}\end{smallmatrix}\bigr)\in (K_m \cap S)s^k$, namely $k=h$.
Therefore
\[
\begin{split}
& \frac{1}{\chi_{ij}(I)}\tau_l(\chi_{ij})\\
=& p^{l+2m}\mu_{\alpha}(\det s^h)e(\mathrm{Tr}s^h)\psi_{ij}(s^h)\\
& \cdot \Big(\sum_{0\leqslant t_1, t_2 <p^{n-m}}\mu_{\alpha}((1+p^mt_1)(1+p^mt_2))\lambda (p^{2m}(i_1t_1+i_2t_2))e(p^mx_h(t_1+t_2))\Big), 
\end{split}
\]
which proves the assertion of the theorem. 
\hfill
$\Box$\\

\subsect{3.5. $\mathfrak{X}_4$}

In this subsection, we consider the Gauss sum $\tau_l(\chi)$ for every irreducible character $\chi$ in $\mathfrak{X}_4$.  
Let $\nu_0$ be an injective multiplicative character such that $\nu_0(1+p^{l-1})=\zeta_p$, and let $\nu=\nu_0\circ\det$. 
By the natural homomorphism $G_l \to G_{l-1}$, we regard $\widehat{G_{l-1}}$ as a subset of $\widehat{G_l}$. 
The derived subgroup of $G_l$ is $SL_2(\Z_{p^l})$, and $G_l/SL_{2}(\Z_{p^l})$ is isomorphic to $\Z_{p^l}^{\times}$. 
Thus the group $\mathrm{Hom}(G_l,\C^{\times})$ is generated by $\nu$. 
Since $\nu^p \in \widehat{G_{l-1}}$,  
the irreducible characters of $G_l$ in $\mathfrak{X}_4$ are $\nu^i \theta \ (0\leqslant i <p,\ \theta \in \widehat{G_{l-1}})$.

Let us evaluate $\tau_l(\nu^i \theta)$ for $0\leqslant i <p$ and $\theta \in \widehat{G_{l-1}}$.  
\[
\begin{split}
& \tau_l(\nu^i \theta) \\
=& \sum_{\begin{smallmatrix}0\leqslant a,b,c,d <p^{l-1},\\ p\nmid ad-bc\end{smallmatrix}}
\sum_{0\leqslant x,y,z,w <p}(\nu^i \theta ) \begin{pmatrix}a & b \\ c& d \end{pmatrix}
\nu \left( I+p^{l-1} \begin{pmatrix}x & y \\ z & w \end{pmatrix} \right)^i \\
&\cdot e\left( \mathrm{Tr}\begin{pmatrix}a & b \\ c& d \end{pmatrix}\left( I+p^{l-1} \begin{pmatrix}x & y \\ z & w \end{pmatrix} \right) \right) \\
=& \sum_{\begin{smallmatrix}0\leqslant a,b,c,d <p^{l-1},\\ p\nmid ad-bc\end{smallmatrix}}
(\nu^i \theta ) \begin{pmatrix}a & b \\ c& d \end{pmatrix}e\left( \mathrm{Tr}\begin{pmatrix}a & b \\ c& d \end{pmatrix}\right)\\
& \cdot \left( \sum_{0\leqslant x < p}\zeta_{p}^{ix}e(p^{l-1}ax) \right)\left( \sum_{0\leqslant w < p}\zeta_{p}^{iw}e(p^{l-1}dw) \right)\\
& \cdot \left( \sum_{0\leqslant y < p}e(p^{l-1}cy) \right)\left( \sum_{0\leqslant z < p}e(p^{l-1}bz) \right)\\
=& p^4 \sum_{\begin{smallmatrix}0\leqslant a,b,c,d <p^{l-1},\\ p\nmid ad-bc, br\equiv cr \equiv 0\ (\mathrm{mod}\ p), \\
i+ar\equiv i+dr \equiv 0\ (\mathrm{mod}\ p)\end{smallmatrix}}
(\nu^i \theta ) \begin{pmatrix}a & b \\ c& d \end{pmatrix}e\left( \mathrm{Tr}\begin{pmatrix}a & b \\ c& d \end{pmatrix}\right) .
\end{split}
\]
Thus, $\tau_l(\nu^i \theta)=0$ if either (i) $1\leqslant i <p$ and $p|r$ or (ii) $i=0$ and $p\nmid r$. \\
If $i=0$ and $p|r$, then we have  
\[
\tau_l(\theta)=p^4\tau_{l-1}(\theta , e),
\]
where $e$ is regarded as an additive character of $\Z_{p^{l-1}}$ on the right hand side.\\
If $1\leqslant i <p$ and $p\nmid r$, then
\[
\tau_l(\nu^i\theta)=p^4 \sum_{\begin{smallmatrix}0\leqslant a,b,c,d <p^{l-1}, b\equiv c \equiv 0\ (\mathrm{mod}\ p), \\
a\equiv d \equiv -ir^{-1}\ (\mathrm{mod}\ p)\end{smallmatrix}}(\nu^i\theta)\begin{pmatrix}a & b \\ c& d \end{pmatrix}e(a+d).
\] 
It seems likely that we cannot simplify this sum any further. 

\sect{4. Appendix}

Let $i$, $j$, $k$ be three integers such that $1\leqslant i \leqslant m$, $1\leqslant j \leqslant i$ and $0\leqslant k \leqslant i$, 
$\lambda$ an injective additive character of $\Z_{p^i}$, and $\beta, b \in \Z_{p^i}^{\times}$. 
Leigh-Cliff-Wen considered the following sum $P$ to determine the values of characters in $\mathfrak{X}_3$ in section $6.1$ of \cite{lcw} p.1314: 
\[
P=\sum_{\begin{smallmatrix}0\leqslant c,d <p^{i},\\ p\nmid d\end{smallmatrix}}\lambda (p^j\beta d+d^{-1}(p^kb-c^2)).
\]
However their evaluation of this sum seems to be mistaken.
The purpose of this section is to give the correct value of this sum. 

Since $\lambda$ is an additive character, we have
\[
P=\sum_{\begin{smallmatrix}0\leqslant d <p^{i},\\ p\nmid d\end{smallmatrix}}\lambda(p^j\beta d+ p^kbd^{-1})\sum_{0\leqslant c<p^i}\lambda(-d^{-1}c^2).
\]
The second sum $\sum_{0\leqslant c<p^i}\lambda(-d^{-1}c^2)$ is the quadratic Gauss sum. 
By the standard evaluation of this sum (cf. \cite{bew} Theorem 1.5.2), we have
\[
\sum_{0\leqslant c<p^i}\lambda(-d^{-1}c^2)=\left(\frac{-rd}{p}\right)^i \left(\frac{-1}{p}\right)^{\delta_i/2}p^{i/2},
\]
where 
$\delta_i=(1-(-1)^i)/2=\begin{cases}1, & \mathrm{if}\ i\ \mathrm{is\ odd}, \\ 0, & \mathrm{if}\ i\ \mathrm{is\ even}, \end{cases}$ 
and $\lambda(1)=\zeta_{p^i}^{r}$. \\
Hence
\[
P=\left(\frac{r}{p}\right)^i\left(\frac{-1}{p}\right)^{i+\delta_i/2}p^{i/2}
\sum_{\begin{smallmatrix}0\leqslant d <p^{i},\\ p\nmid d\end{smallmatrix}}\left(\frac{d}{p}\right)^i\lambda(p^j\beta d+ p^kbd^{-1}).
\]
Put $P_1=\sum_{0\leqslant d <p^{i}, p\nmid d}\left(\frac{d}{p}\right)^i\lambda(p^j\beta d+ p^kbd^{-1})$. 
In order to evaluate $P_1$, we consider the following five cases: \\

\noindent
(i)$j=k=i$\\

In this case, we have
\[
P_1=\sum_{0\leqslant d <p^{i}, p\nmid d}\left(\frac{d}{p}\right)^i
=
\begin{cases}
p^{i-1}(p-1), & \textrm{if $i$ is even}, \\
0, & \textrm{if $i$ is odd}. 
\end{cases} 
\]

\noindent
(ii)$j<k\leqslant i$\\

Let $\lambda_1$ be an additive character of $\Z_{p^{i-j}}$ defined by $\lambda_1(1)=\lambda(p^j)$. Then
\[
\begin{split}
P_1=& p^j\sum_{0\leqslant d <p^{i-j}, p\nmid d}\left(\frac{d}{p}\right)^i\lambda_1(\beta d+ p^{k-j}bd^{-1})\\
=& p^j\left(\frac{\beta}{p}\right)^i\sum_{0\leqslant d <p^{i-j}, p\nmid d}\left(\frac{d}{p}\right)^i\lambda_1(d+ p^{k-j}\beta bd^{-1}).
\end{split}
\]
The map $d\mapsto d+ p^{k-j}\beta bd^{-1}$ is a bijection on $\Z_{p^{i-j}}^{\times}$ and preserves the squaredness. Hence
\[
\begin{split}
P_1=& p^j\left(\frac{\beta}{p}\right)^i\sum_{0\leqslant d <p^{i-j}, p\nmid d}\left(\frac{d}{p}\right)^i\lambda_1(d)\\
=& 
\begin{cases}
-p^{i-1}, & \textrm{if $j=i-1$ and $i$ is even},\\
\left(\frac{\beta r}{p}\right)\left(\frac{-1}{p}\right)^{1/2}p^{i-1/2}, & \textrm{if $j=i-1$ and $i$ is odd},\\
0, & \textrm{otherwise}. 
\end{cases}
\end{split}
\]

\noindent
(iii)$k<j\leqslant i$\\

As in case (ii), we have
\[
P_1= 
\begin{cases}
-p^{i-1}, & \textrm{if $k=i-1$ and $i$ is even},\\
\left(\frac{br}{p}\right)\left(\frac{-1}{p}\right)^{1/2}p^{i-1/2}, & \textrm{if $k=i-1$ and $i$ is odd},\\
0, & \textrm{otherwise}. 
\end{cases}
\]

\noindent
(iv)$j=k<i$ and $i$ is even\\

Let $\lambda_1$ be an additive character defined by $\lambda_1(1)=\lambda(p^j)$. Then 
\[
P_1=p^j\sum_{0\leqslant d <p^{i-j}, p\nmid d}\lambda_1(\beta d+ bd^{-1})\\
=p^jK, 
\]
where $K=\sum_{0\leqslant d <p^{i-j}, p\nmid d}\lambda_1(d+ \beta bd^{-1})$ is the Kloosterman sum. 
While we cannot give an explicit determination of $K$ in the case $i-j=1$,  
there is an evaluation of $K$ in the case $i-j>1$ as follows (cf. \cite{sal} $\S 2$): \\
If $\beta b$ is nonsquared, $K=0$. If $\beta b$ is squared, then
\[
K= \left(\frac{ur}{p}\right)^{j}\left(\frac{-1}{p}\right)^{\delta_j/2}p^{(i-j)/2}\left(\lambda_1(2u)+\left(\frac{-1}{p}\right)^j\lambda_1(-2u)\right),
\]
where $u^2=\beta b$. \\

\noindent
(v)$j=k<i$ and $i$ is odd\\

Let $h=i-j$ and $\lambda_1$ be an additive character of $\Z_{p^{h}}$ defined by $\lambda_1(1)=\lambda(p^j)$. Then
\[
\begin{split}
P_1=& p^j\sum_{0\leqslant d <p^{h}, p\nmid d}\left(\frac{d}{p}\right)\lambda_1(\beta d+ bd^{-1})\\
=& p^j\left(\frac{\beta}{p}\right)\sum_{0\leqslant d <p^{h}, p\nmid d}\left(\frac{d}{p}\right)\lambda_1 ( d+ \beta bd^{-1}).
\end{split}
\]
Put $K^{\prime}=\sum_{0\leqslant d <p^{h}, p\nmid d}\left(\frac{d}{p}\right)\lambda_1 ( d+ \beta bd^{-1})$. 
Since $K^{\prime}=\left(\frac{\beta b}{p}\right)K^{\prime}$, $K^{\prime}=0$ if $\beta b$ is nonsquared. 

Assume $\beta b$ is squared, say $u^2=\beta b$. Let $f$ be a $\C$-valued function on $\Z_{p^h}$ defined by
\[
f(a)=\sum_{0\leqslant d < p^h, p\nmid d}\left(\frac{d}{p}\right)\lambda_1( d+ a^2 d^{-1}).
\]
For $c \in \Z$, we have
\[
\begin{split}
p^h\langle f, \lambda_1^c \rangle_{\Z_{p^h}}= & \sum_{0\leqslant a <p^h}\sum_{0\leqslant d < p^h, p\nmid d}\left(\frac{d}{p}\right)\lambda_1 ( d+ a^2 d^{-1})\lambda_1(-ca)\\
=& \sum_{0\leqslant d < p^h, p\nmid d}\left(\frac{d}{p}\right)\lambda_1(d)\sum_{0\leqslant a <p^h}\lambda_1(d^{-1}(a^2-cda))\\
=& \sum_{0\leqslant d < p^h, p\nmid d}\left(\frac{d}{p}\right)\lambda_1(d)\sum_{0\leqslant a <p^h}\lambda_1(d^{-1}((a-2^{-1}cd)^{2}-(2^{-1}cd)^2)\\
=& \sum_{0\leqslant d < p^h, p\nmid d}\left(\frac{d}{p}\right)\lambda_1 ((1-4^{-1}c^2)d)\sum_{0\leqslant a <p^h}\lambda_1(d^{-1}a^2)\\
=& \left(\frac{r}{p}\right)^{j-1}\left(\frac{-1}{p}\right)^{\delta_{j-1}/2}p^{h/2}\sum_{0\leqslant d <p^h, p\nmid d}\left(\frac{d}{p}\right)^{j}\lambda_1((4-c^2)d).
\end{split}
\]
Put $g(c)=\sum_{0\leqslant d <p^h, p\nmid d}\left(\frac{d}{p}\right)^{j}\lambda_1((4-c^2)d)$. Then
\[
g(c)=
\begin{cases}
p^{h-1}(p-1), & \textrm{if $p^h | 4-c^2$ and $j$ is even},\\
-p^{h-1}, & \textrm{if $4-c^2=p^{h-1}c_1$ for some $p\nmid c_1$ and $j$ is even}, \\
\left(\frac{c_1r}{p}\right)\left(\frac{-1}{p}\right)^{1/2}p^{h-1/2}, & \textrm{if $4-c^2=p^{h-1}c_1$ for some $p\nmid c_1$ and $j$ is odd},\\
0, & \textrm{otherwise}.
\end{cases}
\] 
Hence, if $j$ is even, we have
\[
\begin{split}
& K^{\prime}= f(u)= \sum_{0\leqslant c <p^h}\langle f, \lambda_1^c \rangle_{\Z_{p^h}}\lambda_1(cu)\\
=& \left(\frac{r}{p}\right)\left(\frac{-1}{p}\right)^{1/2}p^{-h/2}\Big( p^{h-1}(p-1)( \lambda_1(2u)+\lambda_1(-2u) )-p^{h-1}
\sum_{\begin{smallmatrix}0\leqslant c <p^h,\\ p^{h-1} \| 4-c^2\end{smallmatrix}}\lambda_1(cu) \Big)\\
=& \left(\frac{r}{p}\right)\left(\frac{-1}{p}\right)^{1/2}p^{(i-j)/2}(\lambda_1 (2u)+\lambda_1 (-2u)),
\end{split}
\]
and if $j$ is odd, we have
\[
\begin{split}
& K^{\prime}= f(u)= \sum_{0\leqslant c <p^h}\langle f, \lambda_1^c \rangle_{\Z_{p^h}}\lambda_1(cu)\\
=& \left(\frac{r}{p}\right)\left(\frac{-1}{p}\right)^{1/2}p^{(h-1)/2}\sum_{1\leqslant e <p}\left( \left(\frac{-e}{p}\right)\lambda_1 ( u(2+p^{h-1}e))+\left(\frac{e}{p}\right)\lambda_1 ( u(-2+p^{h-1}e))\right)\\
=& \left(\frac{r}{p}\right)\left(\frac{-1}{p}\right)^{1/2}p^{(h-1)/2} \left( \left(\frac{-u}{p}\right)\lambda_1(2u)+\left(\frac{u}{p}\right) \lambda_1(-2u)\right) \sum_{1\leqslant e <p} \left(\frac{e}{p}\right)\lambda_1 (p^{h-1}e)\\
=& \left(\frac{u}{p}\right)p^{(i-j)/2}\left(\lambda_1 (2u)+\left(\frac{-1}{p}\right)\lambda_1 (-2u)\right).
\end{split}
\]

\medskip
\medskip
This completes the evaluation of $P_1$ in all cases. 


\end{document}